\documentclass[12pt,oneside,notitlepage]{amsart}
\usepackage{graphicx}
\newtheorem{thm}{Theorem}[section]
\newtheorem{cor}[thm]{Corollary}
\newtheorem{lem}[thm]{Lemma}
\newtheorem{prop}[thm]{Proposition}
\theoremstyle{definition}

\theoremstyle{remark}
\newtheorem{rem}[thm]{Remark}
\numberwithin{equation}{section}

\begin{document}

\title[Stability]{NUMERICAL ALGORITHM FOR FINDING BALANCED METRICS ON VECTOR BUNDLES}%
\author[Seyyedali]{Reza Seyyedali}

\address{Department of Mathematics \\
   Johns Hopkins University \\
   Baltimore, MD 21218}

\email{seyyedali@math.jhu.edu}

\thanks{}%

\begin{abstract}

In \cite{D3}, Donaldson defines a dynamical system on the space of
Fubini-Study metrics on a polarized compact K\"ahler manifold.
Sano proved that if there exists a balanced metric for the
polarization, then this dynamical system always converges to the
balanced metric (\cite{S}).
 In \cite{DKLR}, Douglas, et.\ al., conjecture that the same holds in the case of vector bundles. In this paper, we give an affirmative answer to their conjecture.

\end{abstract}
\maketitle
\pagestyle{plain} \pagenumbering{arabic}

\newtheorem*{Mthm}{Main Theorem}
\newtheorem{Thm}{Theorem}[section]
\newtheorem{Prop}[Thm]{Proposition}
\newtheorem{Lem}[Thm]{Lemma}
\newtheorem{Cor}[Thm]{Corollary}
\newtheorem{Def}[Thm]{Definition}
\newtheorem{Guess}[Thm]{Conjecture}
\newtheorem{Ex}[Thm]{Example}
\newtheorem{Rmk}{Remark}
\newtheorem{Not}{Notation}
\def\thesection{\arabic{section}}
\renewcommand{\theThm} {\thesection.\arabic{Thm}}

\section{Introduction}
In \cite{D3}, Donaldson defines a dynamical system on the space
of Fubini-Study metrics on a polarized compact K\"ahler manifold. Sano proved that if there exists a
balanced metric for the polarization, then this dynamical
system always converges to the balanced metric (\cite{S}).
 In \cite{DKLR}, Douglas, et.\ al., conjecture that the same holds in the case of vector bundles. In this paper, we give an affirmative answer to their conjecture.

Let $(X,\omega)$ be a K\"ahler manifold of dimension $n$ and $E$
be a very ample holomorphic vector bundle on $X$. Let $h$ be a
Hermitian metric on $E$. We can define a $L^2$-inner product on
$H^{0}(X,E)$ by

$$\langle s,t\rangle= \int_{X} h(s,t)\frac{\omega^n}{n!}.$$
Let $s_{1},...,s_{N}$ be an orthonormal basis for $H^{0}(X,E)$
with respect to this $L^2$-inner product. The Bergman kernel of
$h$ is defined by $$B(h)=\sum s_{i}\otimes s_{i}^{*_{h}}.$$ Note
that $B(h)$ does not depend on the choice of the orthonormal basis
$s_{1},...,s_{N}$.

A metric $h$ is called balanced if $B(h)$ is a constant multiple
of the identity. By the theorem of Wang (\cite[Theorem 1.1]{W}),
we know that the existence of balanced metrics is closely related
to the stability of the vector bundle $E$. Indeed $E$ admits a
unique (up to a positive constant) balanced metric if and only if
the Gieseker point of $E$ is stable. On the other hand, a balanced
metric is unique (up to a constant) provided the bundle is simple
(cf. Lemma \ref{newlem2} below). Let $K$ and $M$ be  the space of
Hermitian metrics on $E$ and Hermitian inner product on
$H^{0}(X,E)$ respectively. Following Donaldson (\cite{D2}), one
can define the following maps

\begin{itemize}

\item Define $$\textrm{Hilb}: K \rightarrow M$$ by

 $$\langle s,t \rangle _{\textrm{Hilb}(h)}= \frac{N}{Vr} \int
\langle s(x),t(x) \rangle_{h}\frac{\omega^n}{n!},$$where
$N=\dim(H^{0}(X,E)) $ and $V=\textrm{  Vol}(X,\omega)$. Note that
$\textrm{Hilb}$ only depends on the volume form
$\frac{\omega^{n}}{n!}$.

\item For the metric $H$ in $M_{E}$, $FS(H)$ is the unique metric
on $E$ such that  $\sum s_{i}\otimes s_{i}^{*_{FS(H)}}=I$, where
$s_{1},...,s_{N}$ is an orthonormal basis for $H^{0}(X,E)$ with
respect to $H$. This gives the map $FS:M \rightarrow K$.

\item Define a map $$T:M \rightarrow M$$ by $T(H)=\textrm{Hilb}
\circ FS(H).$ Notice that this map $T$ is called generalized
$T$-operator in \cite{DKLR}.

\end{itemize}

It is easy to see that a metric $h$ is balanced if and only if
$\textrm{Hilb}(h)$ is a fixed point of the map $T$.

The main theorem of this paper is the following
\begin{thm}\label{firstthm}

Suppose that $E$ is simple and admits a balanced metric. Then for
any $H_{0}\in M_{E}$, the sequence $T^r(H_{0})$ converges to
$H_{\infty}$, where $H_{\infty}$ is a balanced metric on $E$.

\end{thm}

Our proof  follows Sano's argument in \cite{S} with the
necessary modifications for the bundle case.

In order to prove the theorem, we consider the functional $Z$ that
is used by Wang (\cite{W}) and Phong, Sturm (\cite{PS}) in order
to study the existence and uniqueness of balanced metrics on
holomorphic vector bundles. The key property of this functional is
that its critical points are balanced metrics. In the first
section we recall some properties of the functionals $Z$ and
$\widetilde{Z}$. In the second section, we give an appropriate
notation of boundedness for subsets of $M$ which is defined in
\cite{S}. It is easy to see that any bounded sequence has a
convergent subsequence after a suitable rescaling of the sequence.
Therefore in order to prove that the sequence $H_{n}=T^n(H)$
converges , we need to show that $H_{n}$ is bounded. On the other
hand, existence of a balanced metric implies that $\widetilde{Z}$
is bounded from below and proper in a suitable sense. Hence it
shows that $\widetilde{Z}(H_{n})$ is bounded. Now properness of
$\widetilde{Z}$ implies that $H_{n}$ is bounded.

\thanks{\textbf{Acknowledgements:} I am sincerely grateful to Richard Wentworth for introducing me the subject and many helpful discussions and suggestions
on the subject. I would also like to thank him for all his help,
support and encouragement.}

\section{Balanced metrics on vector bundles}

As before, let $(X,\omega)$ be a K\"ahler manifold and $E$ be a
very ample holomorphic vector bundle on $X$. Using global sections
of $E$, we can map $X$ into $G(r, H^{0}(X,E)^*)$. Indeed, for any
$x \in X$, we have the evaluation map $H^{0}(X,E)\rightarrow
E_{x}$, which sends $s$ to $s(x)$. Since $E$ is globally
generated, this map is a surjection. So its dual is an inclusion
of $E_{x}^* \hookrightarrow H^{0}(X,E)^*$, which determines a
$r$-dimensional subspace of $H^{0}(X,E)^*$. Therefore we get an
embedding $i: X\hookrightarrow G(r, H^{0}(X,E)^*)$. Clearly we
have $i^*U_{r}=E^*$, where $U_{r}$ is the tautological vector
bundle on $G(r, H^{0}(X,E)^*)$, i.e. at any $r$-plane in $G(r,
H^{0}(X,E)^*)$, the fibre of $U_{r}$ is exactly that $r$-plane. A
choice of basis for $H^{0}(X,E)$ gives an isomorphism between
$G(r, H^{0}(X,E)^*)$ and the standard $G(r,N)$, where $N= \dim
H^{0}(X,E)$. We have the standard Fubini-Study hermitian metric on
$U_{r}$, so we can pull it back to $E$ and get a hermitian metric
on $E$.  Using $i^*h_{FS}$ and $\omega$, we get an $L^{2}$ inner
product on $H^{0}(X,E)$. The embedding is called balanced if
$\int_{X} \, \langle s_{i}, s_{j}\rangle \,
\frac{\omega^{n}}{n!}=C \delta_{ij} $. We can formulate this
definition in terms of maps $\textrm{Hilb}$ and $FS$.

 \begin{Def}

A $\omega$ balanced metric on $E$ is a pair $(h^*,H^*)$ so that
$$\textrm{Hilb}(h^*)=H^*  \textrm{    ,    }   FS(H^*)=h^*$$

 \end{Def}

 Fixing a nonzero element $\Theta \in \bigwedge^N H^{0}(X,E)$,
 We can define the determinant of any element in $M$. Thus we
 can define a map $$\log \det: M\rightarrow \mathbb{R}$$
 A different choice of $\Theta$ only changes this map by an additive
 constant.
 Also, we define a functional $I:K\rightarrow \mathbb{R}$ again
 unique up to an additive constant. Fix a background metric
 $h_{0}$ and consider a path $h_{t}=e^{\phi_{t}}h_{0}$ in $K$
 then \begin{equation}\label{eq1}\frac{dI}{dt}=\int_{X} tr(\dot{\phi}) \textrm{   dVol}_{\omega}\end{equation}
 This functional is a part of Donaldson's functional.
 We define: \begin{equation}\label{eq2} Z=-I \circ FS:M\rightarrow \mathbb{R}\end{equation}

 We have the following scaling identities:
 $$\textrm{Hilb}(e^{\alpha}h)=e^{\alpha} \textrm{Hilb}(h),$$
 $$ FS(e^{\alpha}h)=e^{\alpha} FS(h),$$
 $$I(e^{\alpha}h)=I(h)+\alpha rV  ,$$
 where $\alpha$ is a real number.

Following Donaldson, we define:
\begin{equation}\label{eq3}\widetilde{Z}=Z+\frac{rV}{N}\log\det.\end{equation} So
$\widetilde{Z}$ is invariant under constant scaling of the metric.

This functional $Z$ is studied by Wang in \cite{W} and Phong and
Sturm in \cite{PS}. They consider this as a functional on
$SL(N)/SU(N)$. In order to see this, we observe that there is a
correspondence between $M$ and $GL(N)/U(N)$. Let fix an element
$H_{0} \in M$ and an orthonormal basis $s_{1},...,s_{N}$ for
$H^{0}(X,E)$ with respect to $H_{0}$. Now for any $H \in M$ we
assign $[H(s_{i},s_{j})] \in GL(N)$. Notice that change of the
orthonormal basis only changes this matrix by multiplication by
elements of $U(N)$. So we get a well define element of
$GL(N)/U(N)$. The subset $$M_{0}=\{H \in M | \det [H(s_{i},
s_{j})]=1\}$$ corresponds to $SL(N)/SU(N)$.

We recall the definition of the Gieseker point of the bundle $E$.
We have a natural map
$$T(E): \bigwedge^{ r}H^{0}(X,E)\rightarrow H^{0}(X,\det(E))$$ which
for any $s_{1},..., s_{r}$ in $H^{0}(X,E)$ is defined by
$$T(E)(s_{1}\wedge...\wedge s_{r})(x)=s_{1}(x)\wedge...\wedge s_{r}(x).$$

Since $E$ is globally generated, $T(E)$ is surjective. We can view
$T(E)$ as an element of
$\textrm{Hom}(\bigwedge^{r}H^{0}(X,E),H^{0}(X,\det(E)))$. This is
called the Gieseker point of $E$. Notice that fixing a basis for
$H^0(X,E)$ gives an isomorphism between $\bigwedge^{r}H^{0}(X,E)$
and $\bigwedge^r\mathbb{C}^N$. Hence, there is a natural action of
$GL(N)$ on
$\textrm{Hom}(\bigwedge^{r}H^{0}(X,E),H^{0}(X,\det(E)))$.
Phong-Sturm (\cite{PS}) and Wang (\cite{W}) prove that $Z$ is
convex along geodesics of $SL(N)/SU(N)$ and its critical points
are corresponding to balanced metrics on $E$. Phong and Sturm
prove the following

\begin{thm}(\cite [Theorem 2]{PS})\label{secondthm}
There exists a $SU(N)$- invariant norm $||.||$ on
$\textrm{Hom}(\bigwedge^{r}H^{0}(X,E),H^{0}(X,\det(E)))$ such that
for any $\sigma \in SL(N)$

$$Z(\sigma)= \log \frac{||\sigma .
T(E)||^2}{||T(E)||^2}$$

\end{thm}

\begin{rem}

In \cite{W}, Wang proves a slightly weaker version of Theorem
\ref{secondthm}. He proves that for any norm $||.||$ on
$\textrm{Hom}(\bigwedge^{r}H^{0}(X,E),H^{0}(X,\det(E)))$, there
exists positive constants $c$ and $c'$ such that

$$Z(\sigma) \geq c\log ||\sigma .
T(E)||^2 +c'. $$

\end{rem}

\begin{thm}(\cite[Lemma 3.5]{W},  \cite[Lemma 2.2]{PS})\label{thirdthm}
The functional $Z$ is convex along geodesics of $M$.

\end{thm}

The Kempf-Ness theorem (\cite{KN}) shows that $Z$ is proper and
bounded from below if $T(E)$ is stable under the action of
$SL(N)$.

The following is an immediate consequence of the above theorem and
the fact that balanced metrics are critical points of $Z$. Also
notice that $\widetilde{Z}$ is invariant under the scaling of a
metric by a positive real number.

\begin{thm}\label{fourththm}

Assume that $H_{0}$ is a balanced metric on $E$. Then
$\widetilde{Z}|_{M_{0}}$ is proper and bounded from below.
Moreover $\widetilde{Z}(H)\geq \widetilde{Z}(H_{0})$ for any $H
\in M$.

\end{thm}

\begin{lem}\label{newlemma1}
 For any $H \in M$, we have $$Tr(T(H)H^{-1})=N$$

\end{lem}

\begin{proof}
We define $h=FS(H)$. Let $s_{1},..., s_{N}$ be an $H$-orthonormal
basis. We have,
$$\sum s_{i}\otimes s_{i}^{*_{h}}=I$$
Therefore, $$r=Tr \big( \sum s_{i}\otimes s_{i}^{*_{h}} \big)=
\sum |s_{i}|^2_{h}.$$ Integrating the above equation concludes the
lemma.

\end{proof}

\begin{lem}\label{firstlemma}

For any $H \in M$,
\begin{itemize}

\item $Z(H)\geq Z(T(H)).$
 \item $\log\det (H) \geq\log\det(T(H)).$
\item $\widetilde{Z}(H)\geq \widetilde{Z}(T(H)).$
\end{itemize}

\end{lem}

\begin{proof}

Put $h=FS(H)$ , $H'=\textrm{Hilb}\circ FS(H)$ and
$h'=FS(H')=e^{\varphi}h$. Let $s_{1},..., s_{N}$ be an
$H'$-orthonormal basis.  We have,
$$\sum s_{i}\otimes s_{i}^{*_{h}}=e^{-\varphi}.$$ Hence,
$$\int_{X} tr(-\varphi)=\int_{X} \log\det (e^{-\varphi})\leq
\int_{X} \log \Big(\frac{tr(e^{-\varphi})}{r} \Big)^{r}$$

$$=r \int_{X} \log (tr(e^{-\varphi}))-rV\log r \leq r V\log \Big(\frac{1}{V}  \int_{X} tr(e^{-\varphi})  \Big)-r V\log r $$

$$=r V\log \Big(  \frac{1}{V}\int_{X} \sum |s_{i}|_{h}^{2}  \Big)-r V\log r =0\,\,\,\,\,\,\,\,\,\,$$
This shows the first inequality. For the second one, Lemma
\ref{newlemma1} implies that $tr(H'H^{-1})=N$. Using the
arithmetic -geometric mean inequality, we get
$$\det (H'H^{-1})^{\frac{1}{N}}\leq \frac{tr(H'H^{-1})}{N}=1.$$
This implies that  $\log \det (H'H^{-1}) \leq 0$. The third
inequality is obtained by summing up the first two.


\end{proof}

A bundle $E$ is called simple if $Aut(E)\simeq \mathbb{C}^*$. We
will also need the following

\begin{lem}\label{newlem2}

Suppose that $E$ is simple and admits a balanced metric. Then the
balanced metric is unique up to a positive constant.

\end{lem}

\begin{proof}

Let $H_{\infty}$ be a balanced metric on $E$ and $s_{1},...,
s_{N}$ be an orthonormal basis of $H^{0}(X,E)$ with respect to
$H_{\infty}$. This basis gives an embedding $\iota : X \rightarrow
Gr(r,N)$ such that $\iota^* U_{r}=E$, where $U_{r}\rightarrow
Gr(r,N)$ is the universal bundle over the Grassmannian. Let assume
that $H$ is another element of $M_{0}$. Therefore, there exists an
element $a \in su(N)$ such that $e^{ia} . H_{\infty}=H$. The one
parameter family $\{ e^{ita} \}$ gives a one parameter family of
automorphism of $(Gr(r,N),U_{r})$ and therefore gives a one
parameter family in $\textrm{Aut}(X,E)$. From lemma 3.5 in
\cite{W}, we have
\begin{equation}\label{neq1} \frac{d^2}{dt^2}Z(e^{ita})=\int_{\iota(X)} ||\tilde{a}||^2
dvol_{X},\end{equation} where $\tilde{a}$ is the vector field on
$Gr(r,N)$ generated by the infinitesimal action of $a$ and
$||\tilde{a}||$ is the Fubini-Study metric on $Gr(r,N)$. Suppose
that $H$ is a balanced metric. Therefore it is a minimum for the
functional $Z$. This implies that
$$\frac{d^2}{dt^2}Z(e^{ita})=0.$$
Hence \eqref{neq1} implies that $\tilde{a}\equiv 0.$ This implis
that the one parameter family $\{ e^{ita} \}$ fixes $X$ and
therefore it is a one parameter family of endomorphisms of $E$.
Thus, simplicity of $E$ implies that $ e^{ita}= I_{E}$ which
basically means that $a=0$ and therefore $H=H_{\infty}$.

\end{proof}

\section{Proof of Theorem \ref{firstthm}}

In this section, we follow Sano's argument in (\cite[Section
3]{S}) very closely. Let $s_{1},...,s_{N} $ be a basis for
$H^0(E)$. Using this basis, we can view elements of $M$ as $N
\times N$ matrices. Now using this identification, we state the
following definition stated in Sano (\cite{S}).

\begin{Def}\label{1def} A subset $U \subseteq M$ is called bounded if there exists a number
$R>1$, satisfying the following:

For any $H \in U$, there exists a positive number $\gamma_{H}$ so
that\begin{equation}\label{eq4}\frac{\gamma_{H}}{R} \leq \min
\frac{|H(\xi)|}{|\xi|}\leq \max \frac{|H(\xi)|}{|\xi|} \leq
\gamma_{H} R
\end{equation}

\end{Def}

Note that boundedness does not depend on the choice of the basis.
Also notice that $\min \frac{|H(\xi)|}{|\xi|}$ is the smallest
eigenvalue of the matrix $[H(s_{i},s_{j})]$ and $\max
\frac{|H(\xi)|}{|\xi|}$ is the largest eigenvalue of the matrix
$[H(s_{i},s_{j})]$.




From the definition, one can see that $U$ is bounded if and only
if there exists $R>1$ satisfying the following:

For any $H \in U$, there exists a positive number $\gamma_{H}$ so
that
$$||[H(s_{i},s_{j})]||_{op}\leq \gamma_{H} R     \textrm{   and   }     ||[H(s_{i},s_{j})]^{-1}||_{op}\leq \gamma_{H}^{-1} R.$$

\begin{prop}\label{firstprop}

Any bounded sequence $H_{i}$ has a subsequence $H_{n_{i}}$ such
that $\gamma_{n_{i}}^{-1}H_{n_{i}}$ converges to some point in
$M$. Here $\gamma_{i}=\gamma_{H_{i}}$ in definition \ref{1def}.

\end{prop}

\begin{proof}

The sequence $\gamma_{n_{i}}^{-1}H_{n_{i}}$ is a bounded sequence
in the space of $N \times N$ matrices with respect to the standard
topology. Hence the proposition follows from the fact that the
closure of bounded sets are compact.

\end{proof}

Notice that the standard topology in the space of $N \times N$
matrices is induced by the standard Euclidean norm. Since all
norms on a finite dimensional vector space are equivalent, we can
use the operator norm on the space of $N \times N$ matrices.
Therefore a sequence $\{H_{\alpha}\}$ in $M$ converges to $H \in
M$ if and only if
$$\big |[H_{\alpha}(s_{i},s_{j})]-[H(s_{i},s_{j})]\big |_{op} \rightarrow 0  \,\,\,\, \textrm{    as    }  \,\,\,\, \alpha \rightarrow 0.$$

\begin{lem}\label{secondlemma}

The set $U \subseteq M$ is bounded if and only if there exists a
number $R>1 $ so that for any $H \in U$, we have
$$\frac{1}{R} \leq \min \frac{|\widetilde{H}(\xi)|}{|\xi|}\leq
\max \frac{|\widetilde{H}(\xi)|}{|\xi|} \leq  R,$$ where
$\widetilde{H}=(\det(H))^{-\frac{1}{N}}H.$

\end{lem}

\begin{proof}

Assume that $U$ is bounded. So by definition there exists a number
$R>1$, satisfying the following:

For any $H \in U$, there exists a positive number $\gamma_{H}$ so
that $$\frac{\gamma_{H}}{R} \leq \min \frac{|H(\xi)|}{|\xi|}\leq
\max \frac{|H(\xi)|}{|\xi|} \leq \gamma_{H} R$$ Let $H$ be an
element of $U$. Without loss of generality we can assume that
$H(s_{i},s_{j})=e^{\lambda_{i}}\delta_{ij}$ and $\lambda_{1}\leq
...\leq \lambda_{N}.$ For any $i$, we have $$ \frac{\gamma_{H}}{R}
\leq e^{\lambda_{i}} \leq  \gamma_{H} R.$$ This implies that
$\gamma_{H} \leq R e^{\lambda_{i}}$ and $\gamma_{H} \geq R^{-1}
e^{\lambda_{i}}.$ Therefore
$$e^{\lambda_{N}}\leq \gamma_{H} R \leq R^2 e^{\lambda_{i}},$$ and
$$e^{\lambda_{1}}\geq \gamma_{H} R^{-1} \geq R^{-2} e^{\lambda_{i}},$$
for any $1 \leq i \leq N.$
Hence $$(\det(H))^{\frac{-1}{N}}
e^{\lambda_{N}}= e^{\lambda_{N}-\frac{\sum \lambda_{i}}{N}}=
\Big(\prod e^{\lambda_{N}-\lambda_{i}} \Big)^{\frac{1}{N}} \leq
R^2.
$$
and $$(\det(H))^{\frac{-1}{N}} e^{\lambda_{1}}=
e^{\lambda_{1}-\frac{\sum \lambda_{i}}{N}}= \Big(\prod
e^{\lambda_{1}-\lambda_{i}} \Big)^{\frac{1}{N}} \geq R^{-2}.
$$

\end{proof}

Let $H_{0}$ be an element in $M$. We define the sequence
$\{H_{n}\}$ by $H_{n}=\textrm{Hilb} \circ FS(H_{n-1})$.

\begin{lem}\label{thirdlemma}

If $\{H_{n}\}$ is a bounded sequence in $M$, then $\det(H_{n})$ is
bounded and
$$\det(H_{n+1}H_{n}^{-1})\rightarrow 1 \text{    as         } n\rightarrow \infty.$$

\end{lem}

\begin{proof}
$\widetilde{Z}(H_{n})$ is bounded since the sequence $\{H_{n}\}$
is bounded. On the other hand, lemma \ref{firstlemma} implies that
the sequences $Z(H_{n})$ and $\log\det(H_{n})$ are decreasing. So,
$\log\det(H_{n})$ is bounded and decreasing. Hence, it converges
to some real number. This implies that
$\det(H_{n+1}H_{n}^{-1})\rightarrow 1$ as $n\rightarrow \infty$.

\end{proof}

\begin{lem}\label{fourthlemma}
Assume $\{H_{n}\}$ is a bounded sequence in $M$. Let $H$ be a
fixed element of $M$ and $s_{1}^{(l)},..s_{N}^{(l)}$ be an
orthonormal basis with respect to $H_{l}$ so that the matrix
$[H(s_{i}^{(l)},s_{j}^{(l)})]$ is diagonal.   Then
$$\frac{N}{Vr}\int_{X}\big|s_{i}^{l}\big|_{h_{l}}^2 \,
dvol_{X}\rightarrow 1  \,\,\,\, \textrm{    as }  \,\,\,\,
l\rightarrow\infty,$$ where $h_{n}=FS(H_{n}).$
\end{lem}

\begin{proof}

Let $\hat{s}_{1}^{(l)},..\hat{s}_{N}^{(l)}$ be an orthonormal
basis with respect to $H_{l}$ so that $H_{l+1}(\hat{s}_{i}^{(l)},
\hat{s}_{j}^{(l)} )$ is diagonal. Hence $$\det
\big[H_{l+1}(\hat{s}_{i}^{(l)}, \hat{s}_{j}^{(l)}
)\big]=\prod_{i=1}^{N}H_{l+1}(\hat{s}_{i}^{(l)}, \hat{s}_{1}^{(l)}
).$$ Lemma \ref{thirdlemma} implies that $$\det
\big[H_{l+1}(\hat{s}_{i}^{(l)}, \hat{s}_{j}^{(l)} )\big]
\rightarrow 1.$$ On the other hand lemma \ref{firstlemma} implies
that $$tr \big[ H_{l+1}(\hat{s}_{i}^{(l)},
\hat{s}_{j}^{(l)})\big]=N.$$ We define
$A_{l}(i)=H_{l+1}(\hat{s}_{i}^{(l)}, \hat{s}_{i}^{(l)}).$
Therefore, we have
\begin{equation}\label{eq5}\prod_{i=1}^{N}A_{l}(i)\rightarrow 1 \,\,\,\, \textrm{    as    }  \,\,\,\, l\rightarrow \infty , \end{equation}
\begin{equation}\label{eq6}\sum_{i=1}^{N} A_{l}(i)=N , \,\,\,\, \textrm{  for any } \, 1\leq l \leq N.\end{equation}
We claim that for any
$i$,\begin{equation}\label{eq7}A_{l}(i)\rightarrow 1 \,\,\,\,
\textrm{    as    }  \,\,\,\, l \rightarrow \infty.\end{equation}
Suppose that for some $1 \leq \alpha \leq N$, $\{A_{l}(\alpha) \}
$ does not converge to $1$ as $l \rightarrow \infty.$ This means
that there exists a positive number $\epsilon >0$ and a
subsequence $\{A_{l_{q}}(\alpha)\}$ such that
\begin{equation}\label{eq8}\big|A_{l_{q}}(\alpha)-1 \big| \geq \epsilon.\end{equation}
On the other hand, \eqref{eq6} implies that $A_{l}(i) \leq N$
since $A_{l}(i) \geq 0$ and therefore the sequences
$\{A_{l_{q}}(i)\}$ are bounded for any $1\leq i \leq N$. Hence
there exist nonnegative numbers $A(1),...A(N)$ and a subsequence
$\{ l_{q_{j}} \}$ so that
\begin{equation}\label{eq9}A_{l_{q_{j}}}(i) \rightarrow A(i) \,\,\,\, \textrm{    as    }  \,\,\,\, j\rightarrow \infty.\end{equation}
Therefore, \eqref{eq5}, \eqref{eq6} and \eqref{eq9} imply that
$$\prod_{i=1}^{N } A(i)=1  \,\,\,\,  \textrm{   and   }  \,\,\,\,   \sum _{i=1}^{N } A(i)=N.$$
By arithmetic-geometric mean inequality, we always have
$$\big(\prod_{i=1}^{N} A(i) \big)^{\frac{1}{N}} \leq  \frac{1}{N} \sum _{i=1}^{N } A(i)$$
and equality holds if and only if all $A_{i}$'s are equal. Since
equality holds in this case, we conclude that $A(1)=...=A(N)=1.$
In particular
$$A_{l_{q_{j}}}(\alpha) \rightarrow 1    \,\,\,\, \textrm{    as    }  \,\,\,\,
j\rightarrow \infty,$$ which contradicts \eqref{eq8}. This implies
that $H_{l+1}(\hat{s}_{i}^{(l)}, \hat{s}_{i}^{(l)} )\rightarrow 1$
for all $i$.

On the other hand, there exists $[a_{ij}^l]\in U(N)$ such that
$s_{i}^{(l)}=\sum_{j=1}^{N} a_{ij} \hat{s}_{j}^{l}$. Since $U(N)$
is compact, we can find a subsequence of $[a_{ij}^l]$ which
converges to an element of $U(N)$. Without loss of generality, we
can assume that there exists $[a_{ij}] \in U(N)$ such that
$a_{ij}^{l} \rightarrow a_{ij}$ as $l \rightarrow \infty$. We
have,
$$H_{l+1}(s_{i}^{(l)}, s_{i}^{(l)} )=\sum a_{ij}^l\overline{a_{ik}^l}H_{l+1}(\hat{s}_{j}^{(l)},
\hat{s}_{k}^{(l)} )\rightarrow \sum_{j=1}^{N} |a_{ij}|^2 =1$$

\end{proof}

\begin{prop}(cf. \cite [Proposition ]{S})\label{secondprop}
If $\{H_{n}\}$ is a bounded sequence in $M$, then for any $H \in
M$ and any $\epsilon
>0$,
\begin{equation}\widetilde{Z}(H)> \widetilde{Z}(H_{n})-\epsilon,
\end{equation}
for sufficiently large $n$.

\end{prop}

\begin{proof}

Let  $s_{1}^{(l)},...,s_{N}^{(l)}$ be an orthonormal basis with
respect to $H_{l}$ such that
$H(s_{i}^{(l)},s_{j}^{(l)})=\delta_{ij} e^{\lambda^{(l)}_{i}}$ .
We fix a positive integer $l$.  Define
$H_{t}(s_{i}^{(l)},s_{i}^{(l)})=\delta_{ij}
e^{t\lambda^{(l)}_{i}}$. We have $H_{0}=H_{l}$ and $H_{1}=H$. Let
$f_{l}(t)=f(t)=\widetilde{Z}(H_{t})$. We have
$$f(1)-f(0)=\int_{0}^{1} f'(t)\, dt=\int_{0}^{1}
\Big(f'(0)+\int_{0}^{t} f''(s)\, ds \Big)\, dt$$
$$=f'(0)+\int_{0}^{1} \int_{0}^{1} f''(s) \,ds\,dt\geq
f'(0),$$since $\widetilde{Z}$ is convex along geodesics. On the
other hand, we have
$$f'(t)=\frac{d}{dt}\Big(-I(FS(H_{t}))+\frac{Vr}{N} \log
\det(H_{t})\Big)$$
$$=-\int_{X} \frac{d}{dt}\big(FS(H_{t})\big) dvol_{X}+\frac{Vr}{N}\sum
\lambda^{(l)}_{i}$$ Therefore,
\begin{equation}\label{1eq}f'_{l}(0)=-\int_{X} \big( \sum \lambda^{(l)}_{i}
|s_{i}^{(l)}|^{2}_{h_{l}} \big) dvol_{X}+\frac{Vr}{N}\sum
\lambda^{(l)}_{i},\end{equation} where $h_{l}=FS(H_{l}).$

We have that
$e^{\frac{-\lambda^{(l)}_{1}}{2}}s_{1}^{(l)},...,e^{\frac{-\lambda^{(l)}_{N}}{2}}s_{N}^{(l)}$
is an orthonormal basis with respect to $H$ for any $l$. Hence
lemma \ref{secondlemma} implies that there exists $R >1$ so that
$$\frac{(\det(H_{l}))^{\frac{1}{N}}}{R} <   H_{l}(e^{\frac{-\lambda^{(l)}_{i}}{2}}s_{i}^{(l)},e^{\frac{-\lambda^{(l)}_{i}}{2}}s_{1}^{(i)})            <
(\det(H_{l}))^{\frac{1}{N}}R,$$ for any $i$ and $l$. Therefore
$$ \frac{1}{N} \log (\det(H_{l}))- \log R <-\lambda^{(l)}_{i}<
\frac{1}{N} \log (\det(H_{l}))+ \log R.$$ This implies that
$\{\lambda^{(l)}_{i}\}$ is bounded since $\{\det(H_{l})\}$ is
bounded by Lemma \ref{thirdlemma}. Hence \eqref{1eq} implies that
$f^{'}_{l} (0)\longrightarrow 0$, as $l \longrightarrow \infty.$

\end{proof}

\begin{cor}\label{secondcor}

If $\{H_{n}\}$ is a bounded sequence in $M$, then
$$\widetilde{Z}(H_{n}) \longrightarrow \inf
\{\widetilde{Z}(H)\mid H \in M\}.$$\\

\end{cor}

\begin{proof}[Proof of Theorem \ref{firstthm}]

As before, fix $H_{0} \in M$ and an orthonormal basis
$s_{1},...,s_{N}$ for $H^0(X,E)$ with respect to the metric
$H_{0}$. As in Section $2$, let $$M_{0}=\big \{ H \in M  | \det
[H(s_{i},s_{j})]=1 \big \}.$$ Assume that there exists a balanced
metric on $E$. Since the balanced metric is unique up to a
positive constant, there exists a unique balanced metric
$H_{\infty} \in M_{0}.$ As before, for  any $H \in M$, we define
$$ \widetilde{H}= (\det H)^{-\frac{1}{N}}H.$$ Clearly
$\widetilde{H} \in M_{0}$ and
$$\widetilde{Z}(\widetilde{H})=\widetilde{Z}(H)=Z(\widetilde{H}).$$
Since there exists a balanced metric on $E$, theorem
\ref{fourththm} implies that the functional $Z_{| M_{0}}$ is
proper and bounded from below. Hence the sequence
$Z(\widetilde{H_{n}})$ is a bounded sequence in $\mathbb{R}$ since
the sequence $\widetilde{Z}(H_{n})=Z(\widetilde{H_{n}})$ is
decreasing. Therefore the sequence $\{\widetilde{H_{n}}\}$ is
bounded in $M_{0}$ since $Z_{| M_{0}}$ is proper. We claim that
$$\widetilde{H_{n}}\longrightarrow H_{\infty}   \,\,\,\, \textrm{
as    }  \,\,\,\, n\rightarrow \infty. $$ Suppose that the
sequence $\{\widetilde{H_{n}}\}$ does not converge to
$H_{\infty}$. Then there exists $\epsilon > 0$ and a subsequence
$\{H_{n_{j}}\}$ such that
\begin{equation}\label{eq10}||\widetilde{H_{n_{j}}}- H_{\infty}||_{op}\geq \epsilon.\end{equation}
On the other hand, we know that the sequence
$\{\widetilde{H_{n_{j}}}\}$ is bounded. Therefore there exist a
subsequence $\{\widetilde{H_{n_{j_{q}}}}\}$  and an element
$\widehat{H} \in M $ such that
$$\widetilde{H_{n_{j_{q}}}} \rightarrow \widehat{H}\,\,\,\, \textrm{    as    }  \,\,\,\,    q  \rightarrow \infty.$$
Therefore,
$$1=\det[\widetilde{H_{n_{j_{q}}}}(s_{\alpha}, s_{\beta})] \rightarrow \det[\widehat{H}(s_{\alpha}, s_{\beta})]\,\,\,\, \textrm{    as    }  \,\,\,\, q  \rightarrow \infty,$$
which implies that $\widehat{H} \in M_{0}$.  Now, corollary
\ref{secondcor} implies that
$$\widetilde{Z}(\widetilde{H_{n_{j_{q}}}})=\widetilde{Z}(H_{n_{j_{q}}})
\longrightarrow \inf \{\widetilde{Z}(H)\mid H \in M\}. $$ Hence,
$$\widetilde{Z}(\widehat{H})=\inf \{\widetilde{Z}(H)\mid H \in
M\}.$$ This implies that $\widehat{H}$ is a balanced metric and
therefore $H_{\infty}=\widehat{H}$ by lemma \ref{newlem2}. This
contradicts \eqref{eq10}. Thus $\widetilde{H_{n}}\longrightarrow
H_{\infty} \,\,\,\, \textrm{    as    }  \,\,\,\, q\rightarrow
\infty$.

Now lemma \ref{thirdlemma} implies that $\log \det(H_{n})$ is
bounded. The sequence $\{\log \det(H_{n})\}$ is bounded and
decreasing. Therefore there exists $b \in \mathbb{R}$ such that
$$\log \det(H_{n}) \rightarrow b  \,\,\,\, \textrm{    as    }  \,\,\,\,  n \rightarrow \infty.$$ Hence $\det(H_{n})$ converges to the positive real number $e^b$. Thus $$H_{n}\longrightarrow
e^{\frac{-b}{N}} H_{\infty}    \,\,\,\, \textrm{    as    }
\,\,\,\, n\rightarrow \infty.$$

\end{proof}

\end{document}